\documentclass[a4paper,10pt]{article}
\usepackage{stmaryrd}
\usepackage{amsfonts}
\usepackage{bbm}
\usepackage{amscd}
\usepackage{mathrsfs}
\usepackage{latexsym,amssymb,amsmath,amscd,amscd,amsthm,amsxtra,xypic}
\usepackage[dvips]{graphicx}
\usepackage[utf8]{inputenc}
\usepackage[T1]{fontenc}
\usepackage{lmodern}
\usepackage{amssymb}
\usepackage[all]{xy}
\usepackage{nicefrac,mathtools,enumitem}
\usepackage{microtype}

\newtheorem{thm}{Theorem}[section]
\newtheorem{lem}[thm]{Lemma}
\newtheorem{cor}[thm]{Corollary}
\newtheorem{pro}[thm]{Proposition}
\newtheorem{ex}[thm]{Example}
\newtheorem{rmk}[thm]{Remark}
\newtheorem{defi}[thm]{Definition}

\newcommand {\emptycomment}[1]{} 

\newcommand{\be }{\begin{equation}}
\newcommand{\ee }{\end{equation}}

\newcommand{\pf}{\noindent{\bf Proof.}\ }

\newcommand{\g}{\mathfrak g}

\newcommand{\huaV}{\mathcal{V}}

\newcommand{\frka}{\mathfrak a}
\newcommand{\frkb}{\mathfrak b}
\newcommand{\frkc}{\mathfrak c}
\newcommand{\frkd}{\mathfrak d}
\newcommand{\frke}{\mathfrak e}
\newcommand{\frkf}{\mathfrak f}
\newcommand{\frkg}{\mathfrak g}
\newcommand{\frkh}{\mathfrak h}

\newcommand{\frkl}{\mathfrak l}

\newcommand{\frkX}{\mathfrak X}
\newcommand{\frkY}{\mathfrak Y}

\def\qed{\hfill ~\vrule height6pt width6pt depth0pt}

\newcommand{\Id}{\rm{Id}}
\newcommand{\Fu}{\mathrm{F}}

\newcommand{\br}[1]{   [ \cdot,    \cdot  ]   }

\newcommand{\dM}{\mathrm{d}}

\newcommand{\Der}{\mathrm{Der}}

\newcommand{\gl}{\mathfrak {gl}}

\newcommand{\Ker}{\mathrm{Ker}}

\newcommand{\ad}{\mathrm{ad}}

\newcommand{\ve}{\mathrm{v}}
\newcommand{\Vect}{\mathrm{Vect}}
\newcommand{\sgn}{\mathrm{sgn}}
\newcommand{\Ksgn}{\mathrm{Ksgn}}

\newcommand{\V}{\mathbb{V}}

\begin{document}
\title{
{ $3$-$Lie_\infty$-algebras and $3$-Lie 2-algebras
\thanks
 {
This research is supported by NSFC (11471139) and NSF of Jilin Province (20140520054JH).
 }
} }
\author{Yanqiu Zhou, Yumeng Li and Yunhe Sheng  \\
Department of Mathematics, Jilin University,\\
 Changchun 130012,  China
\\
email:yqzhou15@mails.jlu.edu.cn,~liym1009@mails.jlu.edu.cn,~
shengyh@jlu.edu.cn}
\date{}
\footnotetext{{\it{Keyword}:    $3$-$Lie_\infty$-algebras, $Lod_\infty$-algebras, $3$-Lie
$2$-algebras,  $3$-pre-Lie algebras}} \footnotetext{{\it{MSC}}: 17B99, 53D17.}
\maketitle

\begin{abstract}
In this paper,  we introduce the notions of a  $3$-$Lie_\infty$-algebra and a 3-Lie 2-algebra. The former is a model for a 3-Lie algebra that satisfy the fundamental identity up to all higher homotopies, and the latter is the categorification of a 3-Lie algebra. We prove that the 2-category of 2-term  $3$-$Lie_\infty$-algebras is equivalent to the 2-category of 3-Lie 2-algebras. Skeletal and strict 3-Lie 2-algebras are studied in detail. A construction of a 3-Lie 2-algebra from a symplectic 3-Lie algebra is given.
\end{abstract}


\section{Introduction}
The notion of a Filippov algebra, or an $n$-Lie algebra was introduced in  \cite{Filippov}. It is the algebraic structure corresponding to Nambu mechanics \cite{Gautheron,N,T}. Recently,
due to applications in the Bagger–Lambert–Gustavsson theory of multiple
M2-branes \cite{BL3,BL4,BL2,HHM,P,BL1}, $n$-Lie algebras, or more generally, $n$-Leibniz algebras, are widely studied  \cite{BaiRGuo,CasasPBW,CasasSN,CasasCMn-Lie,CasasHigher,deformation,Kasymov}. See the review article \cite{review} for more details.

Recently, people study higher categorical
structures with motivations from string theory. One way to provide higher categorical
structures is by categorifying existing mathematical concepts. One
of the simplest higher structure is a $2$-vector space, which is a
categorified  vector space. If we further put Lie algebra structures
on $2$-vector spaces, then we obtain the notion of Lie $2$-algebras
\cite{baez:2algebras}. The Jacobi identity is replaced by a natural
transformation, called Jacobiator, which also satisfies some
coherence laws of its own. The 2-category of Lie 2-algebras is equivalent to the 2-category of $L_\infty$-algebras. $L_\infty$-algebras, also called strong homotopy Lie algebras, were introduced  in \cite{Stasheff1}. See \cite{Fialowski,Penkava} for the cohomology and deformation theory of $L_\infty$-algebras. As a  model for ``Leibniz
algebras that satisfy Jacobi identity up to all higher homotopies'',
  the notion of a strongly homotopy (sh)
Leibniz algebra, or a $Lod_\infty$-algebra was given in \cite{livernet} by Livernet,
 which is further studied by Ammar, Poncin and Uchino
in \cite{ammardefiLeibnizalgebra,UchinoshL}.
  In \cite{Leibniz2al}, the authors introduced the notion of a Leibniz 2-algebra, which is the categorification of a Leibniz algebra, and prove that the category of Leibniz 2-algebras and the category of 2-term $Lod_\infty$-algebras are equivalent.

  The aim of this paper is to provide a model for 3-Lie algebras that satisfy the fundamental identity up to all higher homotopies, and give the categorification of 3-Lie algebras. It is well-known that there is a Leibniz algebra structure on the space of fundamental objects associated to a 3-Lie algebra (more generally, $n$-Lie algebras). We define $3$-$Lie_\infty$-algebras in such a way that there is a $Lod_\infty$-algebra structure on the graded vector space of fundamental objects. Furthermore, we put 3-Lie algebra structures on 2-vector spaces and obtain 3-Lie 2-algebras. Similar to the relation between Lie 2-algebras and 2-term $L_\infty$-algebras, we prove that the 2-category of 2-term  $3$-$Lie_\infty$-algebras and the 2-category of 3-Lie 2-algebras are equivalent. Finally, we study skeletal 3-Lie 2-algebras and strict 3-Lie 2-algebras. In particular, we show that strict 3-Lie 2-algebras are equivalent to crossed modules of 3-Lie algebras. The more general notion of crossed modules of $n$-Leibniz had been given in \cite{CasasCMn-Lie}, and the relation with the third cohomology group is established there.

  The paper is organized as follows. In Section 2, we review 3-Lie algebras and their cohomology, 2-vector spaces and $Lod_{\infty}$-algebras. In Section 3, we define $3$-$Lie_\infty$-algebras. We show that given a $3$-$Lie_\infty$-algebra $(\huaV,\{l_{2n+1}\}_{n=0}^\infty)$, we can obtain a $Lod_\infty$-algebra  $(\wedge^2\huaV,\{\frkl_n\}_{n=1}^\infty)$ (Theorem \ref{thm:main}). Then we focus on 2-term $3$-$Lie_\infty$-algebras and show that there is a 2-category of 2-term $3$-$Lie_\infty$-algebras. In Section 4, we give the notion of a 3-Lie 2-algebra, and prove that the 2-category of 2-term  $3$-$Lie_\infty$-algebras is equivalent to the 2-category of 3-Lie 2-algebras (Theorem \ref{thm:equivalence}). In Section 5, first we classify skeletal 3-Lie 2-algebras using the third cohomology group. Then we show that there is a one-to-one correspondence between strict 3-Lie 2-algebras and crossed modules of 3-Lie algebras. Finally, we construct 3-Lie 2-algebras from symplectic 3-Lie algebras via the underlying 3-pre-Lie algebras given in \cite{baiguosheng}.

\section{Preliminaries}

\subsection{3-Lie algebras and their cohomology}

\begin{defi}
A $3$-Lie algebra ${\frkg}$ is a vector space together with a trilinear fully skew-symmetric map $[\cdot,\cdot,\cdot]:\wedge^3 \frkg\longrightarrow \frkg,$ the $3$-bracket, such that the following fundamental identity is satisfied:
\begin{equation}\label{eq:fi}
~[x_1,x_2,[x_3,x_4,x_5]]=[[x_1,x_2,x_3],x_4,x_5]+[x_3,[x_1,x_2,x_4],x_5]+[x_3,x_4,[x_1,x_2,x_5]],\ \forall x_i\in\frkg.
\end{equation}
\end{defi}

  Let $(\g,[\cdot,\cdot,\cdot])$ be a $3$-Lie algebra. For all $x,y\in\g$, define $\ad_{(x,y)}:\g\longrightarrow\g$ by $\ad_{(x,y)}z=[x,y,z]$. Then Eq. \eqref{eq:fi} is equivalent to that $\ad_{(x,y)}$ is a derivation.   Elements in $\wedge^2\g$ are called {\bf fundamental objects} of the $3$-Lie algebra $(\g,[\cdot,\cdot,\cdot])$. There is a bilinear operation $[\cdot,\cdot]_{\Fu}$ on $  \wedge^{2}\g$, which is given by
\begin{equation}\label{eq:bracketfunda}
~[\frkX,\frkY]_{\rm F}=[x_1,x_2,y_1]\wedge y_2+y_1\wedge[x_1,x_2,y_2],\quad \forall \frkX=x_1\wedge x_2, ~\frkY=y_1\wedge y_2.
\end{equation}
It is well-known that $(\wedge^2\g,[\cdot,\cdot]_{\rm F})$ is a Leibniz algebra, which plays important role in the theory of 3-Lie algebras.

\begin{defi}
Let $V$ be a vector space. A representation   of a $3$-Lie algebra $\frkg$ on $V$ is a bilinear map $\rho:\wedge^2\frkg\longrightarrow \gl(V),$ such that
\begin{eqnarray*}[\rho(\frkX),\rho(\frkY)]&=&\rho([\frkX,\frkY]_{\rm F}),\ \forall \frkX,\frkY\in\wedge^2\frkg,\\
 \rho(x,[y_1,y_2,y_3])&=&\rho(y_2,y_3)\rho(x,y_1)-\rho(y_1,y_3)\rho(x,y_2)+\rho(y_1,y_2)\rho(x,y_3),\ \forall x,y_i\in\frkg.
 \end{eqnarray*}
\end{defi}

We denote a representation by $(V;\rho)$. A $p$-cochain on  $\frkg$ with the coefficients in a representation   $(V;\rho)$ is a linear map $\alpha^p:\wedge^2\frkg\otimes\stackrel{(p-1)}{\ldots}\otimes\wedge^2\frkg\wedge\frkg\longrightarrow\frkg$.
Denote the space of $p$-cochains by $C^{p}(\g;V).$ The coboundary operator $\delta:C^{p}(\g;V)\longrightarrow C^{p+1}(\g;V)$  is given by
\begin{eqnarray*}
&&(\delta\alpha^p)(\frkX_1,\ldots,\frkX_p,z)\\
&=& \sum_{1\leq j<k}(-1)^j\alpha^p(\frkX_1,\ldots,\hat{\frkX}_j,\ldots,\frkX_{k-1},[\frkX_j,\frkX_k]_{\rm F},\frkX_{k+1},\ldots,\frkX_{p},z)\\
&&+\sum_{j=1}^p(-1)^j\alpha^p(\frkX_1,\ldots,\hat{\frkX}_j,\ldots,\frkX_{p},[\frkX_j,z])\\
&&+\sum_{j=1}^p(-1)^{j+1}\rho(\frkX_j)\alpha^p(\frkX_1,\ldots,\hat{\frkX}_j,\ldots,\frkX_{p},z)\\
&&+(-1)^{p+1}\Big(\rho(y_{p},z)\alpha^p(\frkX_1,\ldots,\frkX_{p-1},x_{p} ) +\rho(z,x_{p})\alpha^p(\frkX_1,\ldots,\frkX_{p-1},y_{p} ) \Big),
\end{eqnarray*}
for all   $\frkX_i=(x_i,y_i)\in\wedge^2\frkg$ and $z\in\frkg.$

\subsection{2-vector spaces and $Lod_{\infty}$-algebras}\label{2vectorspace and leibnia}
Vector spaces can be categorified to $2$-vector spaces. A good introduction for this subject is \cite{baez:2algebras}.
Let $\Vect$ be the category of vector spaces.

\begin{defi}{\rm\cite{baez:2algebras}}
A $2$-vector space is a category in the category $\Vect$.
\end{defi}

Thus a $2$-vector space $C$ is a category with a vector space of
objects $C_0$ and a vector space of morphisms $C_1$, such that all
the structure maps are linear. Let $s,t:C_1\longrightarrow C_0$ be
the source and target maps respectively. Let $\cdot_\ve$ be the
composition of morphisms. Let $1:C_0\longrightarrow C_1$ be the unit map, i.e. for all $x\in C_0$, $1_x\in C_1$ is the identity morphism from $x$ to $x$.

It is well known that the 2-category of 2-vector spaces is
equivalent to the 2-category of 2-term complexes of vector spaces.
Roughly speaking, given a 2-vector space $C$, we have a 2-term complex of vector spaces
\begin{equation}\label{eq:complex}
\Ker(s)\stackrel{t}{\longrightarrow}C_0.
\end{equation}
Conversely, any 2-term complex of vector spaces
$\huaV:V_1\stackrel{\dM}{\longrightarrow}V_0$ gives rise to a
2-vector space of which the set of objects is $V_0$, the set of
morphisms is $V_0\oplus V_1$, the source map $s$ is given by
$s(v+m)=v$, and the target map $t$ is given by $t(v+m)=v+\dM m$,
where $v\in V_0,~m\in V_1.$ We denote the 2-vector space associated
to the 2-term complex of vector spaces
$\huaV:V_1\stackrel{\dM}{\longrightarrow}V_0$ by $\V$:
\begin{equation}\label{eqn:V}
\V=\begin{array}{c}
\V_1:=V_0\oplus V_1\\
\vcenter{\rlap{s }}~\Big\downarrow\Big\downarrow\vcenter{\rlap{t }}\\
\V_0:=V_0.
 \end{array}\end{equation}

The notion of  $Lod_\infty$-algebras, also called strongly homotopy (sh) Leibniz algebras, was introduced in  \cite{livernet} and well studied in \cite{ammardefiLeibnizalgebra,UchinoshL}.

\begin{defi}{\rm}
A $Lod_\infty$-algebra is a graded  vector space $L=L_0\oplus
L_1\oplus\cdots$ equipped with a system $\{l_k|~1\leq k<\infty\}$ of
linear maps $l_k:\otimes^kL\longrightarrow L$ with degree
$\deg(l_k)=k-2$, where the exterior powers are interpreted in the
graded sense and the following relation  is satisfied:
\begin{eqnarray*}
&&\sum_{i+j=n+1}\sum_{j\leq k\leq n
}\sum_{\sigma}(-1)^{(k+1-j)(j-1)}(-1)^{j(|x_{\sigma(1)}|+\cdots+|x_{\sigma(k-j)}|)}\sum_{\sigma}\sgn(\sigma)\Ksgn(\sigma)\\
&&l_i(x_{\sigma(1)},\dots,x_{\sigma(k-j)},l_j(x_{\sigma(k+1-j)},\dots,x_{\sigma(k-1)},x_k),x_{ k+1},\dots,x_{n})=0,
\end{eqnarray*}
where the summation is taken over all $(k-j,j-1)$-unshuffles and
``$\Ksgn(\sigma)$'' is the Koszul sign for a permutation $\sigma\in
S_k$, i.e. $$ x_1\wedge x_2\wedge\cdots\wedge
x_k=\Ksgn(\sigma)x_{\sigma(1)}\wedge x_{\sigma(2)}\wedge\cdots\wedge
x_{\sigma(k)}.
$$
\end{defi}
We will denote a $Lod_\infty$-algebra by $(\huaV,\{l_{k}\}_{k=1}^\infty)$.

\section{3-$Lie_\infty$-algebras}\label{sec:pre}

In this section, we give the notion of a $3$-$Lie_\infty$-algebra and show that associated to a $3$-$Lie_\infty$-algebra, there is naturally a $Lod_\infty$-algebra structure on the graded space of fundamental objects. Furthermore, we focus on 2-term $3$-$Lie_\infty$-algebras and construct the corresponding 2-category.

For all $n\in \mathbb N$  and $1\leq j\leq k\leq n$, we say that $\sigma\in S_{k}$ is a $(k-j,j)$-unshuffle if
$$\sigma(1)\leq\sigma(2)\leq\cdots\leq \sigma(k-j)\ \ and\ \ \sigma(k-j+1)\leq\sigma(k-j+2)\leq\cdots\leq \sigma(k).$$

\begin{defi}\label{defi:liewuqiong}
A $3$-$Lie_\infty$-algebra is a graded vector space $\huaV=V_0\oplus V_1\oplus\cdots $, equipped with a system $\{l_{2n+1}:n=0,1,2,\cdots\}$ of linear maps: $l_{2n+1}:(\wedge^2 \huaV)\otimes \stackrel{n}{\ldots}\otimes(\wedge^2 \huaV)\wedge\huaV\longrightarrow \huaV,$ with degree $n-1,$  such that for all $X_a=(x_a,y_a)\in \wedge^2 \huaV$ and $x\in \huaV$, the following relation is satisfied:
\begin{eqnarray}
\nonumber&&\sum_{i+j=n-1}\sum_{j\leq k \leq n-1
}\sum_{\sigma\in(k-j-1,j)\mbox{-unshuffle}}(-1)^{(k-j)(j)}(-1)^{(j+1)(|X_{\sigma(1)}|+\cdots+|X_{\sigma(k-j-1)}|)}\sgn(\sigma)\Ksgn(\sigma)\\
\nonumber&&l_{2i+1}(X_{\sigma(1)},\dots,X_{\sigma(k-j-1)},\widetilde{l}_{2j+1}(X_{\sigma(k-j)},\dots,X_{\sigma(k-1)},X_{k}),X_{k+1},\dots,X_{n-1},x)\\
\nonumber&&+\sum_{i+j=n-1}\sum_{\sigma\in (i,j)\mbox{-unshuffle}}(-1)^{(i+1)(j)}(-1)^{(j+1)(|X_{\sigma(1)}|+\cdots+|X_{\sigma(i)}|)}\sgn(\sigma)\Ksgn(\sigma)\\
\label{eq:gFI}&&l_{2i+1}(X_{\sigma(1)},\dots,X_{\sigma(i)},l_{2j+1}(X_{\sigma(i+1)},\dots,X_{\sigma(n-1)},x))=0,
\end{eqnarray}
where  $\widetilde{l}_{2n+1}:\otimes^{n+1}(\wedge^2 \huaV)\longrightarrow \wedge^2 \huaV$ is induced by $l_{2n+1}$ via
\begin{eqnarray}
\nonumber \widetilde{l}_{2n+1}(X_1,X_2,\cdots,X_{n+1})&\triangleq &
\label{eq:formularln} (-1)^{|x_{n+1}|(|X_1|+\ldots+|X_n|+n-1)}x_{n+1}\wedge l_{2n+1}(X_1,\cdots,X_n,y_{n+1})\\
&&+l_{2n+1}(X_{1},\cdots,X_n,x_{n+1})\wedge y_{n+1}.
\end{eqnarray}
\end{defi}

We will denote a 3-$Lie_\infty$-algebra by $(\huaV,\{l_{2n+1}\}_{n=0}^\infty)$. By \eqref{eq:gFI}, we have $l_1\circ l_1=0$, which makes $(\huaV,l_1)$ being a complex of vector spaces. Thus, we will write $\dM=l_1$ sometimes.


Similar as the fact that there is a Leibniz algebra structure on the space of fundamental objects associated to a 3-Lie algebra, we have a $Lod_\infty$-algebra structure on the graded space of fundamental objects associated to a 3-$Lie_\infty$-algebra.

\begin{thm}\label{thm:main}
Let $(\huaV,\{l_{2n+1}\}_{n=0}^\infty)$ be a $3$-$Lie_\infty$-algebra. Then, $(\wedge^2\huaV,\{\frkl_n\}_{n=1}^\infty)$ is a $Lod_\infty$-algebra, where $\frkl_n=\widetilde{l}_{2n-1}$ is given by \eqref{eq:formularln}.
\end{thm}
\pf Denote the left hand side of \eqref{eq:gFI} by $\Xi(X_1,\cdots,X_{n-1},x)$, we have
{\footnotesize
\begin{eqnarray*}
&&\sum_{i+j=n+1}\sum_{1\leq j\leq k\leq n
}\sum_{\sigma}(-1)^{(k+1-j)(j-1)}(-1)^{j(|X_{\sigma(1)}|+\cdots+|X_{\sigma(k-j)}|)}\sum_{\sigma}\sgn(\sigma)\Ksgn(\sigma)\\
&&\frkl_i(X_{\sigma(1)},\dots,X_{\sigma(k-j)},\frkl_j(X_{\sigma(k+1-j)},\dots,X_{\sigma(k-1)},X_{k}),X_{k+1},\dots,X_{n})\\
&=&\sum_{i+j=n+1}\sum_{j\leq k\leq n-1}\sum_{\sigma}(-1)^{(k+1-j)(j-1)}(-1)^{j(|X_{\sigma(1)}|+\cdots+|X_{\sigma(k-j)}|)}\sum_{\sigma}\sgn(\sigma)\Ksgn(\sigma)\\
&&\Big(l_{2i-1}\big(X_{\sigma(1)},\dots,X_{\sigma(k-j)},\widetilde{l}_{2j-1}(X_{\sigma(k+1-j)},\dots,X_{\sigma(k-1)},X_{k}),X_{k+1},\dots,X_{i+j-2},x_{n}\big)\wedge y_{n}\\
&&+(-1)^{(|X_{\sigma(1)}|+\cdots +|X_{\sigma(n-1)}|+i-2+j-2)|x_{n}|} \\
&&x_{n}\wedge l_{2i-1}\big(X_{\sigma(1)},\dots,X_{\sigma(k-j)},\widetilde{l}_{2j-1}(X_{\sigma(k+1-j)},\dots,X_{\sigma(k-1)},X_{k}),X_{k+1},\dots,X_{n-1},y_{n}\big)\Big) \\
&&+\sum_{i+j=n+1} \sum_{\sigma}(-1)^{i(j-1)}(-1)^{j(|X_{\sigma(1)}|+\cdots+|X_{\sigma(i-1)}|)}\sum_{\sigma}\sgn(\sigma)\Ksgn(\sigma)\\
&&\widetilde{l}_{2i-1}\Big(X_{\sigma(1)},\dots,X_{\sigma(i-1)},l_{2j-1}(X_{\sigma(i)},\dots,X_{\sigma(n-1)},x_{n})\wedge y_{n}\\
&&+(-1)^{|x_{n}|(|X_{\sigma(i)}|+\dots+|X_{\sigma(n-1)}|+j-2)}x_{n}\wedge l_{2j-1}(X_{\sigma(i)},\dots,X_{\sigma(n-1)},y_{n})\Big)\\
&=&\sum_{i+j=n-1}\sum_{0\leq j\leq k\leq
n-1}\sum_{\sigma}(-1)^{(k-j)j}(-1)^{(j+1)(|X_{\sigma(1)}|+\cdots+|X_{\sigma(k-j-1)}|)}\sum_{\sigma}\sgn(\sigma)\Ksgn(\sigma)\\
&&l_{2i+1}(X_{\sigma(1)},\dots,X_{\sigma(k-j-1)},\widetilde{l}_{2j+1}(X_{\sigma(k-j)},\dots,X_{\sigma(k-1)},X_{k}),X_{k+1},\dots,X_{n-1},x_{n})\wedge y_{n}\\
&&+(-1)^{(|X_{\sigma(1)}|+\cdots +|X_{\sigma(i+j)}|+i+j-2)|x_{i+j+1}|}\\
&&x_{n}\wedge l_{2i+1}(X_{\sigma(1)},\dots,X_{\sigma(k-j-1)},\widetilde{l}_{2j+1}(X_{\sigma(k-j)},\dots,X_{\sigma(k-1)},X_{k}),X_{k+1},\dots,X_{n-1},y_{n})\\
&&+\sum_{i+j=n-1} \sum_{\sigma}(-1)^{(i+1)j}(-1)^{(j+1)(|X_{\sigma(1)}|+\cdots+|X_{\sigma(i)}|)}\sum_{\sigma}\sgn(\sigma)\Ksgn(\sigma)\\
&&\{l_{2i+1}\Big(X_{\sigma(1)},\dots,X_{\sigma(i)},l_{2j+1}(X_{\sigma(i+1)},\dots,X_{\sigma(n-1)},x_{n})\Big)\wedge y_{n}\\
&&+(-1)^{|x_{n}|(|X_{\sigma(i+1)}|+\dots+|X_{\sigma(n-1)}|+j-1)+|x_{n}|(|X_{\sigma(1)}|+\cdots|X_{\sigma(i)}|+i-1 )}\\
&&x_{n}\wedge l_{2i+1}\Big(X_{\sigma(1)},\dots,X_{\sigma(i)},l_{2j+1}(X_{\sigma(i+1)},\dots,X_{\sigma(n-1)},y_{n})\Big)\}\\
&=&\Xi(X_1,\cdots,X_{n-1},x_{n})\wedge y_{n}\\
&&+(-1)^{(|X_{\sigma(1)}|+\cdots +|X_{\sigma(n-1)}|+i+j-2)|x_{n}|}x_{n}\wedge \Xi(X_1,\cdots,X_{n-1},y_{n})\\
&=&0.
\end{eqnarray*}
}
Therefore, $(\wedge^2\huaV,\frkl_n)$ is a $Lod_\infty$-algebra. \qed\vspace{3mm}

 In particular, if we concentrate on the $2$-term case, we can give explicit formulas for $2$-term $3$-$Lie_\infty$-algebras as follows:

\begin{lem}\label{lem:2term 3L}
A $2$-term $3$-$Lie_\infty$-$algebra$ $\huaV=(V_1,V_0,\dM,l_3,l_5),$ consists of the following data:
\begin{itemize}
\item[$\bullet$] a complex of vector spaces $V_1\stackrel{\dM}{\longrightarrow}V_0,$

\item[$\bullet$] completely skew-symmetric trilinear maps $l_3:V_i\times V_j\times V_k\longrightarrow
V_{i+j+k}$, where $0\leq i+j+k\leq 1$,

\item[$\bullet$] a  multilinear map $l_5:(V_0\wedge V_0)\otimes (V_0\wedge V_0\wedge V_0)\longrightarrow
V_1$,
   \end{itemize}
   such that for any $x,y,x_i\in V_0$ and $f,g,h\in V_1$, the following equalities are satisfied:
\begin{itemize}
\item[$\rm(a)$] $\dM l_3(x,y,f)=l_3(x,y,\dM f),$
\item[$\rm(b)$]$l_3(f,g,h)=0; \quad l_3(f,g,x)=0,$
\item[$\rm(c)$]$l_3(\dM f,g,x)=l_3(f,\dM g,x),$
\item[$\rm(d)$]
$\dM l_5(x_1,x_2,x_3,x_4,x_5)=-l_3(x_1,x_2,l_3(x_3,x_4,x_5))+l_3(x_3,l_3(x_1,x_2,x_4),x_5)\\
\hspace{5cm}\qquad  \qquad \qquad \qquad +l_3(l_3(x_1,x_2,x_3),x_4,x_5)+l_3(x_3,x_4,l_3(x_1,x_2,x_5)),$
\item[$\rm(e)$]$  l_5(\dM f,x_2,x_3,x_4,x_5)=-l_3(f,x_2,l_3(x_3,x_4,x_5))+l_3(x_3,l_3(f,x_2,x_4),x_5)\\
\qquad  \qquad \qquad \qquad+l_3(l_3(f,x_2,x_3),x_4,x_5)+l_3(x_3,x_4,l_3(f,x_2,x_5)),$
\item[$\rm(f)$]$  l_5(x_1,x_2,\dM f,x_4,x_5)=-l_3(x_1,x_2,l_3(f,x_4,x_5))+l_3(f,l_3(x_1,x_2,x_4),x_5)\\
\qquad  \qquad \qquad \qquad+l_3(l_3(x_1,x_2,f),x_4,x_5)+l_3(f,x_4,l_3(x_1,x_2,x_5)),$
\item[$\rm(g)$]
{\footnotesize
\begin{eqnarray*}
&&l_3(l_5(x_1,x_2,x_3,x_4,x_5),x_6,x_7)+l_3(x_5,l_5(x_1,x_2,x_3,x_4,x_6),x_7)+l_3(x_1,x_2,l_5(x_3,x_4,x_5,x_6,x_7))\\
&+&l_3(x_5,x_6,l_5(x_1,x_2,x_3,x_4,x_7))+l_5(x_1,x_2,l_3(x_3,x_4,x_5),x_6,x_7)+l_5(x_1,x_2,x_5,l_3(x_3,x_4,x_6),x_7)\\
&+&l_5(x_1,x_2,x_5,x_6,l_3(x_3,x_4,x_7))=l_3(x_3,x_4,l_5(x_1,x_2,x_5,x_6,x_7))+l_5(l_3(x_1,x_2,x_3),x_4,x_5,x_6,x_7)\\
&+&l_5(x_3,l_3(x_1,x_2,x_4),x_5,x_6,x_7)+l_5(x_3,x_4,l_3(x_1,x_2,x_5),x_6,x_7)+l_5(x_3,x_4,x_5,l_3(x_1,x_2,x_6),x_7)\\
&+&l_5(x_1,x_2,x_3,x_4,l_3(x_5,x_6,x_7))+l_5(x_3,x_4,x_5,x_6,l_3(x_1,x_2,x_7)).
\end{eqnarray*}
}
   \end{itemize}
\end{lem}
Equations (a) and (c) tells us how the differential $\dM$ and the  bracket $l_3$ interact. Equations (d), (e) and (f) tell us that the fundamental identity no longer holds on the nose, but controlled by $l_5$. Equation (g) gives the coherence law that $l_5$ should satisfy.

\begin{cor} Let $(V_1,V_0,\dM,l_3,l_5)$ be a $2$-term $3$-$Lie_\infty$-algebras. Then, for all $f,g,h\in V_1,$ we have
 \begin{equation}\label{eq:3dd}
   l_3(\dM f,\dM g,h)=l_3(\dM f,g,\dM h)=l_3(f,\dM g,\dM h)
 \end{equation}
\end{cor}

We continue by defining homomorphisms between $2$-term $3$-$Lie_\infty$-algebras:

\begin{defi}\label{defi:Lwuqiong hom}
Let $\huaV=(V_1,V_0,\dM,l_3,l_5)$ and $\huaV'=(V_1',V_0',\dM',l_3',l_5')$ be $2$-term $3$-$Lie_\infty$-algebras. A  homomorphism $\phi:\huaV \longrightarrow \huaV'$ consists of:
\begin{itemize}
\item[$\bullet$] a chain map $\phi:\huaV \longrightarrow \huaV'$, which consists of linear maps $\phi_0:V_0 \longrightarrow V_0'$ and $\phi_1:V_1 \longrightarrow V_1'$ preserving the differential;
\item[$\bullet$] a completely skew-symmetric trilinear map $\phi_2:V_0 \times V_0 \times V_0 \longrightarrow V_1'$,
\end{itemize}
such that  for all $x_i\in V_0$ and $h\in V_1,$ we have
 \begin{eqnarray}
 \label{eq:homo1}\dM' (\phi_2(x_1,x_2,x_3))&=&\phi_0(l_3(x_1,x_2,x_3))-l_3'(\phi_0(x_1),\phi_0(x_2),\phi_0(x_3)),\\
 \label{eq:homo2} \phi_2(x_1,x_2,\dM h)&=&\phi_1(l_3(x_1,x_2,h))-l_3'(\phi_0(x_1),\phi_0(x_2),\phi_1(h)),
 \end{eqnarray}
and
    \begin{eqnarray}
    \nonumber&&l_5'(\phi_0(x_1),\phi_0(x_2),\phi_0(x_3),\phi_0(x_4),\phi_0(x_5))+l_3'(\phi_2(x_1,x_2,x_3),\phi_0(x_4),\phi_0(x_5))\\
    \nonumber &+&l_3'(\phi_0(x_3),\phi_2(x_1,x_2,x_4),\phi_0(x_5))+l_3'(\phi_0(x_3),\phi_0(x_4),\phi_2(x_1,x_2,x_5))\\
    \nonumber &+&\phi_2(l_3(x_1,x_2,x_3),x_4,x_5)+\phi_2(x_3,l_3(x_1,x_2,x_4),x_5)+\phi_2(x_3,x_4,l_3(x_1,x_2,x_5))\\
    \label{eq:homo3} &=&l_3'(\phi_0(x_1),\phi_0(x_2),\phi_2(x_3,x_4,x_5))+\phi_2(x_1,x_2,l_3(x_3,x_4,x_5))+\phi_1(l_5(x_1,x_2,x_3,x_4,x_5)).
    \end{eqnarray}
\end{defi}


Let $\varphi:\huaV \longrightarrow \huaV'$ and $\psi:\huaV' \longrightarrow \huaV''$ be $3$-$Lie_\infty$-homomorphisms, their {\bf composition } $((\varphi\circ\psi)_0,(\varphi\circ\psi)_1,(\varphi\circ\psi)_2)$ is given by $(\varphi\circ\psi)_0=\varphi_0\circ\psi_0$,  $(\varphi\circ\psi)_1=\varphi_1\circ\psi_1$, and
 $$(\varphi\circ\psi)_2(x,y,z)=\psi_2(\varphi_0(x),\varphi_0(y),\varphi_0(z))+\psi_1(\varphi_2(x,y,z)).$$
  The {\bf identity homomorphism} $1_{\huaV}:\huaV\longrightarrow\huaV$ has the identity chain map as its underlying map, together with $(1_\huaV)_2=0$.

\begin{defi}\label{defi:Lwuqiong 2hom}
Let $\huaV$ and $\huaV'$ be $2$-term $3$-$Lie_\infty$-algebras, and $\varphi,\psi :\huaV \longrightarrow \huaV'$ be $3$-$Lie_\infty$-homomorphisms. A $3$-$Lie_\infty$-$2$-homomorphism $\tau:\varphi\Rightarrow\psi$ is a chain homotopy such that for all $x_1,x_2,x_3\in V_0 $, the following equation holds:
\begin{eqnarray}\label{2homoinfinity}
\nonumber(\varphi_2-\psi_2)(x_1,x_2,x_3)&=&l_3'(\varphi_0(x_1),\varphi_0(x_2),\tau(x_3)) +l_3'(\dM'\tau(x_1),\tau(x_2), \varphi_0(x_3))+c.p.\\&&-\tau(l_3(x_1,x_2,x_3))
+l_3'(\dM'\tau(x_1),\dM'\tau(x_2),\tau(x_3)).
\end{eqnarray}
\end{defi}

Now we define the vertical and horizontal composition for these $2$-homomorphisms.
Let $\huaV$, $\huaV'$ be $2$-term $3$-$Lie_\infty$-algebras, and $\varphi,\psi,\mu:\huaV\longrightarrow \huaV'$ be $3$-$Lie_\infty$-homomorphisms. Let $\tau:\varphi\Rightarrow \psi$ and $\tau':\psi\Rightarrow \mu$ be $3$-$Lie_\infty$-$2$-homomorphisms.
The {\bf vertical composition} of $\tau$ and $\tau'$, denoted by $\tau'\tau$, is given by $\tau'\tau=\tau'+\tau$.

Let $\huaV$, $\huaV'$, $\huaV''$ be $2$-term $3$-$Lie_\infty$-algebras,   $\varphi,\psi,:\huaV\longrightarrow \huaV'$ and $\varphi',\psi':V'\longrightarrow \huaV''$  $3$-$Lie_\infty$-homomorphisms, and $\tau:\varphi\Rightarrow \psi$ and $\tau':\varphi'\Rightarrow \psi'$  $3$-$Lie_\infty$-$2$-homomorphisms.
The {\bf horizontal composition} of $\tau$ and $\tau'$, denoted by $\tau'\circ\tau$, is given by $\tau'\circ\tau(x)=\tau'_{\varphi_0(x)}+\varphi_1'\tau(x)$.

 Finally, given a $3$-$Lie_\infty$-homomorphism $\varphi$, the {\bf identity $2$-homomorphism} $1_{\varphi}:\varphi\Rightarrow \varphi$ is the zero chain homotopy  $1_{\varphi}(x)=0.$


It is straightforward to see that
\begin{pro}
 There is a  $2$-category {\rm \bf $2$Term$3$-Lie$_\infty$} with $2$-term $3$-$Lie_\infty$-algebras as objects, $3$-$Lie_\infty$-homomorphisms as morphisms, $3$-$Lie_\infty$-$2$-homomorphisms as $2$-morphisms.
\end{pro}

\section{  3-Lie 2-algebras}\label{sec:pre}

In this section, we define 3-Lie 2-algebras, which are the categorification of 3-Lie algebras, and show that the 2-category of 3-Lie 2-algebras is equivalent to the 2-category of 2-term $3$-$Lie_\infty$-algebras.

\begin{defi}\label{defi:3lie2algebra}
A   $3$-Lie $2$-algebra consists of:
\begin{itemize}
\item[$\bullet$] a $2$-vector spaces $L$;
\item[$\bullet$] a completely skew-symmetric trilinear functor, the $\textbf{bracket}$, $[\cdot,\cdot,\cdot]: L\times L\times L\longrightarrow L$;
\item[$\bullet$] a multilinear natural isomorphism $J_{x_1,x_2,x_3,x_4,x_5}$ for all $x_i\in L_0$,
$$[x_1,x_2,[x_3,x_4,x_5]]\stackrel{J_{x_1,x_2,x_3,x_4,x_5}}{\longrightarrow }[[x_1,x_2,x_3],x_4,x_5]+[x_3,[x_1,x_2,x_4],x_5]+[x_3,x_4,[x_1,x_2,x_5]],$$
\end{itemize}
such that for all $x_1,\ldots,x_7\in L_0$, the following fundamentor identity holds:
$$[x_1,x_2,J_{x_3,x_4,x_5,x_6,x_7}](J_{x_1,x_2,[x_3,x_4,x_5],x_6,x_7}+J_{x_1,x_2,x_5,[x_3,x_4,x_6],x_7}+J_{x_1,x_2,x_5,x_6,[x_3,x_4,x_7]})$$
$$([x_5,x_6,J_{x_1,x_2,x_3,x_4,x_7}]+1)([x_5,J_{x_1,x_2,x_3,x_4,x_6},x_7]+[J_{x_1,x_2,x_3,x_4,x_5},x_6,x_7]+1)=$$
$$J_{x_1,x_2,x_3,x_4,[x_5,x_6,x_7]}([x_3,x_4,J_{x_1,x_2,x_5,x_6,x_7}]+1)(J_{[x_1,x_2,x_3],x_4,x_5,x_6,x_7}+J_{x_3,[x_1,x_2,x_4],x_5,x_6,x_7}+$$
\begin{equation}\label{Jacobiator indentity}
J_{x_3,x_4,[x_1,x_2,x_5],x_6,x_7}+J_{x_3,x_4,x_5,[x_1,x_2,x_6],x_7}+J_{x_3,x_4,x_5,x_6,[x_1,x_2,x_7]}),
\end{equation}
or, in terms of a commutative diagram,
$$
\xymatrix{&[x_1,x_2,[x_3,x_4,[x_5,x_6,x_7]]]\ar[dr]^{J_{x_1,x_2,x_3,x_4,[x_5,x_6,x_7]}}\ar[dl]_{[x_1,x_2,J_{x_3,x_4,x_5,x_6,x_7}]}&\\
\frka\ar[d]_{P}&&
\frkc\ar[d]^{1}\\
\frkb\ar[d]_{T}&&C\ar[d]^{S}\\
\frkd\ar[dr]_{Q}&&
\frke\ar[dl]^{R}\\
&\frkf,}
$$
where $\frka,\frkb,\frkc,\frkd,\frke,\frkf$ and $P,Q,R,S,T$ are given by
{\footnotesize
\begin{eqnarray*}
\frka&=&[x_1,x_2,[[x_3,x_4,x_5],x_6,x_7]]+[x_1,x_2,[x_5,[x_3,x_4,x_6],x_7]]+[x_1,x_2,[x_5,x_6,[x_3,x_4,x_7]]],\\
\frkb&=&[[x_1,x_2,[x_3,x_4,x_5]],x_6,x_7]+[[x_3,x_4,x_5],[x_1,x_2,x_6],x_7]+[[x_3,x_4,x_5],x_6,[x_1,x_2,x_7]]\\
&&+[[x_1,x_2,x_5],[x_3,x_4,x_6],x_7]+[x_5,[x_1,x_2,[x_3,x_4,x_6]],x_7]+[x_5,[x_3,x_4,x_6],[x_1,x_2,x_7]]\\
&&+[[x_1,x_2,x_5],x_6,[x_3,x_4,x_7]]+[x_5,[x_1,x_2,x_6],[x_3,x_4,x_7]]+[x_5,x_6,[x_1,x_2,[x_3,x_4,x_7]]],\\
\frkc&=&[x_3,x_4,[x_1,x_2,[x_5,x_6,x_7]]]+[[x_1,x_2,x_3],x_4,[x_5,x_6,x_7]]+[x_3,[x_1,x_2,x_4],[x_5,x_6,x_7]],\\
\frkd&=&[[x_1,x_2,[x_3,x_4,x_5]],x_6,x_7]+[[x_3,x_4,x_5],[x_1,x_2,x_6],x_7]+[[x_3,x_4,x_5],x_6,[x_1,x_2,x_7]]\\
&&+[[x_1,x_2,x_5],[x_3,x_4,x_6],x_7]+[x_5,[x_1,x_2,[x_3,x_4,x_6]],x_7]+[x_5,[x_3,x_4,x_6],[x_1,x_2,x_7]]\\
&&+[[x_1,x_2,x_5],x_6,[x_3,x_4,x_7]]+[x_5,[x_1,x_2,x_6],[x_3,x_4,x_7]]+[x_5,x_6,[[x_1,x_2,x_3],x_4,x_7]]\\
&&+[x_5,x_6,[x_3,[x_1,x_2,x_4],x_7]]+[x_5,x_6,[x_3,x_4,[x_1,x_2,x_7]]],\\
\frke&=&[[x_1,x_2,x_3],x_4,[x_5,x_6,x_7]]+[x_3,[x_1,x_2,x_4],[x_5,x_6,x_7]]+[x_3,x_4,[[x_1,x_2,x_5],x_6,x_7]]\\
&&+[x_3,x_4,[x_5,[x_1,x_2,x_6],x_7]]+[x_3,x_4,[x_5,x_6,[x_1,x_2,x_7]]],\\
\frkf&=&[[[x_1,x_2,x_3],x_4,x_5],x_6,x_7]+[[x_3,[x_1,x_2,x_4],x_5],x_6,x_7]+[[x_3,x_4,[x_1,x_2,x_5]],x_6,x_7]\\
&&+[x_5,[[x_1,x_2,x_3],x_4,x_6],x_7]+[x_5,[x_3,[x_1,x_2,x_4],x_6],x_7]+[x_5,[x_3,x_4,[x_1,x_2,x_6]],x_7]\\
&&+[[x_3,x_4,x_5],[x_1,x_2,x_6],x_7]+[[x_3,x_4,x_5],x_6,[x_1,x_2,x_7]]+[[x_1,x_2,x_5],[x_3,x_4,x_6],x_7]\\
&&+[x_5,[x_3,x_4,x_6],[x_1,x_2,x_7]]+[[x_1,x_2,x_5],x_6,[x_3,x_4,x_7]]+[x_5,[x_1,x_2,x_6],[x_3,x_4,x_7]]\\
&&+[x_5,x_6,[[x_1,x_2,x_3],x_4,x_7]]+[x_5,x_6,[x_3,[x_1,x_2,x_4],x_7]]+[x_5,x_6,[x_3,x_4,[x_1,x_2,x_7]]],\\
P&=&J_{x_1,x_2,[x_3,x_4,x_5],x_6,x_7}+J_{x_1,x_2,x_5,[x_3,x_4,x_6],x_7}+J_{x_1,x_2,x_5,x_6,[x_3,x_4,x_7]},\\
T&=&[x_5,x_6,J_{x_1,x_2,x_3,x_4,x_7}]+1,\\
Q&=&[x_5,J_{x_1,x_2,x_3,x_4,x_6},x_7]+[J_{x_1,x_2,x_3,x_4,x_5},x_6,x_7]+1,\\
S&=&[x_3,x_4,J_{x_1,x_2,x_5,x_6,x_7}]+1,\\
R&=&J_{[x_1,x_2,x_3],x_4,x_5,x_6,x_7}+J_{x_3,[x_1,x_2,x_4],x_5,x_6,x_7}+J_{x_3,x_4,[x_1,x_2,x_5],x_6,x_7}\\
&&+J_{x_3,x_4,x_5,[x_1,x_2,x_6],x_7}+J_{x_3,x_4,x_5,x_6,[x_1,x_2,x_7]}.
\end{eqnarray*}
}
\end{defi}

We continue by setting up a $2$-category of  $3$-Lie $2$-algebras.

\begin{defi}\label{defi:3liehomo}
  Given $3$-Lie $2$-algebras $L$ and $L'$, a homomorphism $F:L\longrightarrow L'$ consists of:
\begin{itemize}
\item[$\bullet$] A linear functor $F$ from the underlying $2$-vector space of $L$ to that of $L',$ and

\item[$\bullet$] a completely skew-symmetric trilinear natural transformation
 $$F_2(x,y,z):[F_0(x),F_0(y),F_0(z)]'\longrightarrow F_0[x,y,z]$$
 such that the following diagram commutes:
$$
\xymatrix{[F_0(x_1),F_0(x_2),[F_0(x_3),F_0(x_4),F_0(x_5)]']'\ar[d]^{[1,1,F_2]}\ar[rrrr]^{\qquad \qquad \qquad \qquad \qquad \qquad \qquad \qquad \qquad J'_{F_0(x_1),F_0(x_2),F_0(x_3),F_0(x_4),F_0(x_5)}\qquad \qquad \qquad \qquad \qquad}&&&&\frka\ar[d]_{[F_2,1,1]+[1,F_2,1]+[1,1,F_2]}\\
[F_0(x_1),F_0(x_2),F_0[x_3,x_4,x_5]]'\ar[d]^{F_2}&&&&
\frkb\ar[d]_{F_2+F_2+F_2}\\
F_0[x_1,x_2,[x_3,x_4,x_5]]\ar[rrrr]^{\qquad F_1(J_{x_1,x_2,x_3,x_4,x_5})}&&&&\frkc\ .}
$$
where $\frka,\frkb,\frkc$ are given by
\begin{eqnarray*}
\frka&=&[[F_0(x_1),F_0(x_2),F_0(x_3)]',F_0(x_4),F_0(x_5)]'+[F_0(x_3),[F_0(x_1),F_0(x_2),F_0(x_4)]',F_0(x_5)]'\\
&&+[F_0(x_3),F_0(x_4),[F_0(x_1),F_0(x_2),F_0(x_5)]']',\\
\frkb&=&[F_0[x_1,x_2,x_3],F_0(x_4),F_0(x_5)]'+[F_0(x_3),F_0[x_1,x_2,x_4],F_0(x_5)]'\\
&&+[F_0(x_3),F_0(x_4),F_0[x_1,x_2,x_5]]',\\
\frkc&=&F_0[[x_1,x_2,x_3],x_4,x_5]+F_0[x_3,[x_1,x_2,x_4],x_5]+F_0[x_3,x_4,[x_1,x_2,x_5]].
\end{eqnarray*}
\end{itemize}
\end{defi}

The identity homomorphism $\Id_L:L\longrightarrow L$ has the identity functor as its underlying functor, together with an identity natural transformation as $(\Id_L)_2.$ Let $L,L'$ and $L''$ be $3$-Lie $2$-algebras, the composite of $3$-Lie $2$-algebra homomorphisms $F:L\longrightarrow L'$ and $G:L'\longrightarrow L''$ which we denote by $G\circ F$, is given by letting the functor $((G\circ F)_0,(G\circ F)_1)$ be the usual composition of $(G_0,G_1)$ and $(F_0,F_1)$, and letting $(G\circ F)_2$ be the following composite:

$$\xymatrix{
  [G_0\circ F_0(x),G_0\circ F_0(y),G_0\circ F_0(z)]'' \ar[dd]_{G_2(F_0(x),F_0(y),F_0(z))} \ar[dr]^{\ \ \ \ (G\circ F)_2(x,y,z)} \\& G_0\circ F_0[x,y,z] .  \\
  G_0[F_0(x),F_0(y),F_0(z)]'  \ar[ur]_{G_1(F_2(x,y,z))}                     }
  $$

We also have $2$-homomorphisms between homomorphisms:

\begin{defi}\label{defi:3lie2homo}

Let $F,G:L\longrightarrow L'$ be $3$-Lie $2$-algebra homomorphisms. A $2$-homomorphism $\theta:F\Rightarrow G$ is a linear natural transformation from $F$ to $G$ such that the following diagram commutes:
$$
\xymatrix{[F_0(x),F_0(y),F_0(z)]'\ar[d]^{[\theta_x,\theta_y,\theta_z]'}\ar[rr]^{F_2}&&F_0[x,y,z]\ar[d]_{\theta_{[x,y,z]}}\\
[G_0(x),G_0(y),G_0(z)]'\ar[rr]^{G_2}&&
G_0[x,y,z].}
$$
\end{defi}

Since $2$-homomorphisms are just natural transformations with an extra property, we vertically and horizontally compose these in the usual way, and an identity $2$-homomorphism is just an identity natural transformation.

It is straightforward to see that

\begin{pro}
There is a   $2$-category   {\rm \bf $3$Lie$2$Alg} with   $3$-Lie $2$-algebras as objects, $3$-Lie $2$-algebra homomorphisms  as morphisms, and $3$-Lie $2$-algebra $2$-homomorphisms as $2$-morphisms.
\end{pro}

Now we  establish the equivalence between the 2-category of $3$-Lie $2$-algebras and that of $2$-term $3$-$Lie_\infty$-algebras. This result is based on the equivalence between $2$-vector spaces and $2$-term chain complexes described in Subsection \ref{2vectorspace and leibnia}.

\begin{thm}\label{thm:equivalence}
The $2$-categories {\rm \bf $2$Term$3$-Lie$_\infty$} and {\rm \bf $3$Lie$2$Alg} are $2$-equivalent.\end{thm}
\pf First we construct a 2-functor $T:$ {\rm \bf $2$Term$3$-Lie$_\infty$} $\longrightarrow$ {\rm \bf $3$Lie$2$Alg}. Given a $2$-term $3$-$Lie_\infty$-algebra $\huaV=(V_1,V_0,\dM,l_2,l_3)$, we have a $2$-vector space $L$ via \eqref{eqn:V}. More precisely, $L_0=V_0, L_1=V_0\oplus V_1$, and the source and the target map are given by $s(x+f)=x$ and $t(x+f)=x+\dM f$. Define
  a skew-symmetric trilinear functor $[\cdot,\cdot,\cdot]:L\times L\times L\longrightarrow L$ by
  \begin{eqnarray*}
  [x+f,y+g,z+h]&=&l_3(x,y,z)+l_3(x,y,h)+l_3(x,g,z)+l_3(f,y,z)\\
  &&+l_3(\dM f,g,z)+l_3(\dM f,y,h)+l_3(x,\dM g,h)+l_3(\dM f,\dM g,h),
  \end{eqnarray*}
and define the fundamentor $J_{x_1,x_2,x_3,x_4,x_5}$ by
$$J_{x_1,x_2,x_3,x_4,x_5}=([x_1,x_2,[x_3,x_4,x_5]],l_5(x_1,x_2,x_3,x_4,x_5)).$$
The source of $J_{x_1,x_2,x_3,x_4,x_5}$ is $[x_1,x_2,[x_3,x_4,x_5]]$ as desired. By (d) in the Lemma \ref{lem:2term 3L}, its target is
\begin{eqnarray*}
t(J_{x_1,x_2,x_3,x_4,x_5})&=&[x_1,x_2,[x_3,x_4,x_5]]+\dM l_5(x_1,x_2,x_3,x_4,x_5)\\
&=&[[x_1,x_2,x_3],x_4,x_5]+[x_3,[x_1,x_2,x_4],x_5]+[x_3,x_4,[x_1,x_2,x_5]],
\end{eqnarray*}
 as desired.
 By Conditions (e) and (f) in Lemma \ref{lem:2term 3L}, we deduce that $J$ is a natural transformation.
By Condition (h) in Lemma \ref{lem:2term 3L}, we can deduce that the fundamentor identity holds.
This completes the construction of a $3$-Lie $2$-algebra $L=T(\huaV)$ from a $2$-term $3$-$Lie_\infty$-algebra $\huaV.$

 We go on to construct a $3$-Lie $2$-algebra homomorphism $T(\phi):T(\huaV)\longrightarrow T(\huaV')$ from a $3$-$Lie_\infty$-algebra homomorphism $\phi=(\phi_0,\phi_1,\phi_2):\huaV\longrightarrow \huaV'$ between $2$-term $3$-$Lie_\infty$-algebras.
Let $T(\huaV)=L$ and $T(\huaV')=L'.$ We define the underlying linear functor of $T(\phi)=F$ with $F_0=\phi_0$, $F_1=\phi_0\oplus\phi_1.$ Define $F_2:V_0\times V_0\times V_0\longrightarrow V_0'\oplus V_1'$ by
$$F_2(x_1,x_2,x_3)=([\phi_0(x_1),\phi_0(x_2),\phi_0(x_3)]',\phi_2(x_1,x_2,x_3)).$$
Then $F_2(x_1,x_2,x_3)$ is a natural isomorphism from $[F_0(x_1),F_0(x_2),F_0(x_3)]'$ to $F_0[x_1,x_2,x_3]$, and $F=(F_0,F_1,F_2)$ is a homomorphism from $L$ to $L'$.
We can also prove that $T$ preserve identities and composition of homomorphisms. So $T$ is a functor.

Furthermore, to construct T to be a 2-functor, we only need to define T on 2-morphisms. Let $\varphi,\psi:\huaV\longrightarrow\huaV'$ be homomorphisms and $\tau:\varphi\Rightarrow\psi$ a 2-homomorphism. Then we define  $$\theta(x)=T(\tau)(x)=(\varphi_0(x),\tau(x)).$$ By \eqref{2homoinfinity}, $T(\tau)$ is a 2-homomorphism from $T(\phi)$ to $T(\psi)$. It is obvious that $T$ preserves the compositions and identities.
\emptycomment{
$$\xymatrix{
  x \ar[d]_{f} \ar[r]^{\theta}
                & \theta(x) \ar[d]^{F_1}  \\
  x' \ar[r]_{\theta}
                & \theta(x')             }$$
                }
Thus, $T$ is a 2-functor from {\rm \bf $2$Term$3$-Lie$_\infty$} to {\rm \bf $3$Lie$2$Alg}.

Next we   construct a 2-functor $S:${\rm \bf $3$Lie$2$Alg} $\longrightarrow$ {\rm \bf $2$Term$3$-Lie$_\infty.$}
Given a $3$-Lie $2$-algebra $L$, we obtain a complex of vector spaces  $\huaV=S(L)$ via \eqref{eq:complex}. More precisely, $V_1=\ker(s),V_0=L_0$ and $\dM=t|_{\ker(s)}$.
For all $x_1,x_2,x_3,x_4,x_5\in V_0=L_0$ and $f,g,h\in V_1\subseteq L_1$,
we define  $l_3$ and $l_5$ as follows:
\begin{itemize}
\item[$\rm(1)$]$l_3(x_1,x_2,x_3)=[1_{x_1},1_{x_2},1_{x_3}],$
\item[$\rm(2)$]$l_3(x_1,x_2,h)=[1_{x_1},1_{x_2},h],$
\item[$\rm(3)$]$l_3(x_1,f,h)=0,~l_3(f,h,g)=0,$
\item[$\rm(4)$]$l_5(x_1,x_2,x_3,x_4,x_5)=J_{x_1,x_2,x_3,x_4,x_5}-1_{s(J_{x_1,x_2,x_3,x_4,x_5})}.$
\end{itemize}
The
various conditions of $L$ being a 3-Lie 2-algebra imply that $\huaV=(V_1,V_0,\dM,l_3,l_5)$ is 2-term 3-$Lie_\infty$-algebra.
This completes the construction of a $2$-term $3$-$Lie_\infty$-algebra $\huaV=S(L)$ from a $3$-Lie $2$-algebra $L$.

Let $L$ and $L'$ be 3-Lie 2-algebras, and $F=(F_0,F_1,F_2):L\longrightarrow L'$ a homomorphism. Let $S(L)=\huaV$ and $S(L')=\huaV'$. We go on to  construct a $3$-$Lie_\infty$-homomorphism $\phi=S(F):\huaV\longrightarrow \huaV'$.
Let $\phi_0=F_0$, $\phi_1=F_1|_{\ker(s)}$. Define $\phi_2:V_0\times V_0\times V_0\longrightarrow V'_1$ by
$$\phi_2(x_1,x_2,x_3)=F_2(x_1,x_2,x_3)-1_{s(F_2(x_1,x_2,x_3))}.$$
Then $\phi_2$ is  completely skew-symmetric, and
 \begin{eqnarray*}
\dM' \phi_2(x_1,x_2,x_3)&=&(t'-s')F_2(x_1,x_2,x_3)\\
&=&\phi_0(l_3(x_1,x_2,x_3))-l_3'(\phi_0(x_1),\phi_0(x_2),\phi_0(x_3)).
\end{eqnarray*}
 The naturality of $F_2$ gives Equation \eqref{eq:homo2} in Definition \ref{defi:Lwuqiong hom}, and the fundamentor identity gives Equation \eqref{eq:homo3} in Definition \ref{defi:Lwuqiong hom}. Thus, $\phi=S(F)$ is a homomorphism between 2-term $3$-$Lie_\infty$-algebras.

Let $F,G:L\longrightarrow L'$ be $3$-Lie $2$-algebra homomorphisms and  $\theta:F\Rightarrow G$  a $2$-homomorphism.
Let $\varphi=S(F),\psi=S(G) :\huaV \longrightarrow \huaV'$ be the corresponding $3$-$Lie_\infty$-homomorphisms. We define
$$\tau(x)=S(\theta)(x)=\theta(x)-1_{s'(\theta(x))}.$$
 By the commutative diagram in Definition \ref{defi:3lie2homo}, we can deduce that \eqref{2homoinfinity} holds. Thus, $\tau=S(\theta)$  is a 2-homomorphism. It is straightforward to deduce that $S$ preserves the compositions and identities. Thus  $S$ is a 2-functor from {\rm \bf $3$Lie$2$Alg} to {\rm \bf $2$Term$3$-Lie$_\infty$}.

In the end, it is easy to construct the natural isomorphisms $\alpha:ST\Rightarrow 1_{\rm\bf 3Lie2Alg}$ and $\beta:TS\Rightarrow 1_{\rm\bf 2Term3Lie_\infty}$. We omit details. The proof is completed. \qed

\section{Skeletal and strict 3-Lie 2-algebras}\label{sec:pre}

By Theorem \ref{thm:equivalence}, we see that 3-Lie 2-algebras and 2-term 3-$Lie_\infty$-algebras are equivalent. Thus, we will call a 2-term 3-$Lie_\infty$-algebra a 3-Lie 2-algebra in the sequel.

A 3-Lie 2-algebra $(V_1,V_0,\dM,l_3,l_5)$ is called {\bf skeletal} ({\bf strict})  if $\dM=0$ ($l_5=0$).

In this section, first we classify skeletal 3-Lie 2-algebras via the third cohomology group. Then, we introduce the notion of a crossed module of 3-Lie algebras, and show that they are equivalent to strict 3-Lie 2-algebra. At the end, we construct a 3-Lie 2-algebra from a 3-Lie algebra with a symplectic structure.

\begin{thm}
There is a one-to-one correspondence between skeletal $3$-Lie $2$-algebras  and quadruples $((\g,[\cdot,\cdot,\cdot]),V,\rho,\Theta)$, where $(\g,[\cdot,\cdot,\cdot])$ is a $3$-Lie algebra, $V$ is a vector space, $\rho$ is a representation of $\g$ on $V,$ and $\Theta$ is a $3$-cocycle on $\g$ with values in $V.$
\end{thm}
\pf Let $(V_1,V_0,\dM=0,l_3,l_5)$ be a skeletal 3-Lie 2-algebras. By (d) in Lemma \ref{lem:2term 3L}, we see that $l_3|_{V_0}$ satisfies the fundamental identity. Thus, $(V_0,l_3|_{V_0})$ is a $3$-Lie algebra.  $l_3$ also gives rise to a map $\rho:\wedge^2V_0\longrightarrow V_1$ by
 \begin{equation}
 \rho(x_1, x_2)(f)=l_3(x_1,x_2,f).
 \end{equation}
   By (e) and (f) in Lemma \ref{lem:2term 3L}, we deduce that $\rho$ is a representation of the 3-Lie algebra $(V_0,l_3|_{V_0})$ on $V_1$. Finally, by (g) in Lemma \ref{lem:2term 3L}, we get that $l_{5}$ is a 3-cocycle.

  Conversely, given a $3$-Lie algebra $(\g,[\cdot,\cdot,\cdot]),$ a representation $\rho$ of $\g$ on a  vector space $V,$ and a $3$-cocycle $\Theta$ on $\g$ with values in $V,$ we define $V_0=\g,$ $V_1=V$,  $\dM=0,$ and totally skew-symmetric $l_3$, $l_5$   by
  $$
  l_3(x,y,z)=[x,y,z],\quad l_3(x,y,f)=\rho(x,y)(f), \quad l_5=\Theta.
  $$
 Then, it is straightforward to see that $(V_1,V_0,\dM=0,l_3,l_5)$ is a skeletal 3-Lie 2-algebra. \qed\vspace{3mm}

Now we introduce the notion of a crossed module of 3-Lie algebras and show that they are equivalent to strict 3-Lie 2-algebras.
\begin{defi}
A crossed module of $3$-Lie algebras is a quadruple  $((\frkg,[\cdot,\cdot,\cdot]_{\frkg}),(\frkh,[\cdot,\cdot,\cdot]_{\frkh}),\mu,\alpha)$, where $(\frkg,[\cdot,\cdot,\cdot]_{\frkg})$ and $(\frkh,[\cdot,\cdot,\cdot]_{\frkh})$ are $3$-Lie algebras, $\mu:\frkg\longrightarrow\frkh$  is a homomorphism of $3$-Lie algebras, and   $\alpha:\wedge^2\frkh\longrightarrow\Der(\frkg)$ is a representation, such that for all $x,y,z\in \frkh,f,g,h\in\frkg,$ the following equalities hold:
\begin{eqnarray}
  \label{eq:cmc1}\mu(\alpha(x,y)(f))&=&[x,y,\mu(f)]_{\frkh},\\
  \label{eq:cmc2}\alpha(\mu(f),\mu(g))(h)&=&[f,g,h]_{\frkg},\\
  \label{eq:cmc3}\alpha(x,\mu(f))(g)&=&-\alpha(x,\mu(g))(f).
\end{eqnarray}
\end{defi}

\begin{rmk}
  The more general notion of crossed modules of $n$-Leibniz had been given in \cite{CasasCMn-Lie}, and the relation with the third cohomology group is established there.
\end{rmk}

\begin{thm}
There is a one-to-one correspondence between strict $3$-Lie $2$-algebras and crossed modules of $3$-Lie algebras.
\end{thm}
\pf  Let $(V_1,V_0,\dM,l_3,l_5)$ be a strict 3-Lie 2-algebras. Define $\frkg=V_1,\frkh=V_0,$ and the following two bracket operations on $\frkg$ and $\frkh$:
\begin{eqnarray}
~[f,g,h]_{\frkg}&=&l_3(\dM f,\dM g,h)=l_3(\dM f,g,\dM h)=l_3(f,\dM g,\dM h),\\
~[x,y,z]_{\frkh}&=&l_3(x,y,z).
\end{eqnarray}
It is straightforward to see that both $[\cdot,\cdot,\cdot]_{\frkg}$ and $[\cdot,\cdot,\cdot]_{\frkh}$ are $3$-Lie brackets. Let $\mu=\dM,$ by Condition (a) in Lemma \ref{lem:2term 3L}, we have
$$\mu[f,g,h]_{\frkg}=\dM l_3(\dM f,\dM g,h)=l_3(\dM f,\dM g,\dM h)=[\mu(f),\mu(g),\mu(h)]_{\frkh},$$
which implies that $\mu$ is a homomorphism of $3$-Lie algebras. Define   $\alpha:\wedge^2\frkh\longrightarrow\gl(\frkg)$ by
$$\alpha(x,y)(f)=l_3(x,y,f).$$
By Condition $(f)$ in Lemma \ref{lem:2term 3L}, we deduce that $\alpha(x,y)\in\Der(\g)$. By Conditions $(e)$ and $(f)$ in Lemma \ref{lem:2term 3L}, we can obtain that $\alpha$ is a representation. Furthermore, we have
\begin{eqnarray*}
  \mu(\alpha(x,y)(f))&=&\dM (\alpha(x,y)(f))=\dM l_3(x,y,f)=l_3(x,y,\dM f)= [x,y,\mu(f)]_{\frkh},\\
  \alpha(\mu(f),\mu(g))(h)&=&l_3(\dM f,\dM g,h)=[f,g,h]_{\frkg},\\
  \alpha(x,\mu(f))(g)&=&l_3(x,\mu(f),g)=-l_3(x,\mu(g),f)=-\alpha(x,\mu(g))(f).
\end{eqnarray*}
Therefore we obtain a crossed module of $3$-Lie algebras.

Conversely, a crossed module of $3$-Lie algebras gives rise to a $2$-term $3$-$Lie_\infty$-algebra, in which $V_1=\frkg,$  $V_0=\frkh$, $\dM=\mu$, and the totally skew-symmetric trilinear map $l_3$ is given by
\begin{eqnarray*}
l_3(x,y,z)&=&[x,y,z]_{\frkh},\\
l_3(x,y,f)&=&\alpha(x,y)(f).
\end{eqnarray*}
The crossed module conditions give various conditions for $2$-term $3$-$Lie_\infty$-algebras. We omit details. The proof is completed. \qed\vspace{3mm}

Now we give an example of crossed module of 3-Lie algebras.

\begin{ex}{\rm
  Let $\g=\mathbb R^3$ with a basis $\{e_1,e_2,e_3\}$, and the $3$-Lie bracket is given by $[e_1,e_2,e_3]_\g=e_1$. Let $\frkh=\mathbb R^3$ with a basis $\{e^1,e^2,e^3\}$, and the $3$-Lie bracket is given by $[e^1,e^2,e^3]_\frkh=e^1$. Define $\mu:\frkg\longrightarrow\frkh$ by
  $\mu(e_i)=e^i$. Obviously, $\mu$ is an homomorphism. Define $\alpha:\wedge^2\frkh\longrightarrow\gl(\g)$ by
  $$
  \alpha(e^1,e^2)=\left(\begin{array}{ccc}
  0&0&1\\
  0&0&0\\
  0&0&0
  \end{array}\right),\quad  \alpha(e^2,e^3)=\left(\begin{array}{ccc}
 1&0&0\\
  0&0&0\\
  0&0&0
  \end{array}\right), \quad  \alpha(e^3,e^1)=\left(\begin{array}{ccc}
  0&1&0\\
  0&0&0\\
  0&0&0
  \end{array}\right).
  $$
  More precisely, $\alpha(e^1,e^2)(e_3)=e_1,~\alpha(e^2,e^3)(e_1)=e_1,~\alpha(e^3,e^1)(e_2)=e_1$. It is straightforward to deduce that $\alpha$ is a representation and take values in $\Der(\g)$. Furthermore, \eqref{eq:cmc1}-\eqref{eq:cmc3} are also satisfied. Thus, $(\g,\frkh,\mu,\alpha)$ constructed above is a crossed module of 3-Lie algebras.
}
\end{ex}

At the end of this section, we construct a strict 3-Lie 2-algebra from a 3-Lie algebra with a symplectic structure.  The main ingredient in this construction is the underlying 3-pre-Lie algebra structure.

\begin{defi}{\rm \cite{baiguosheng}}
 Let $A$ be a vector space with a multilinear map $\{\cdot,\cdot,\cdot\}:A\otimes A\otimes A\rightarrow A$.
$(A,\{\cdot,\cdot,\cdot\})$ is called a {\bf $3$-pre-Lie algebra} if the following identities hold:
\begin{eqnarray}
\{x,y,z\}&=&-\{y,x,z\},\label{eq:d1}\\
\nonumber\{x_1,x_2,\{x_3,x_4,x_5\}\}&=&\{[x_1,x_2,x_3]_C,x_4,x_5\}+\{x_3,[x_1,x_2,x_4]_C,x_5\}\\
&&+\{x_3,x_4,\{x_1,x_2,x_5\}\},\label{eq:d2}\\
\nonumber\{ [x_1,x_2,x_3]_C,x_4, x_5\}&=&\{x_1,x_2,\{ x_3,x_4, x_5\}\}+\{x_2,x_3,\{ x_1,x_4,x_5\}\}\\
&&+\{x_3,x_1,\{ x_2,x_4, x_5\}\}.\label{eq:d3}
\end{eqnarray}
 where $[\cdot,\cdot,\cdot]_C$ is given by
\begin{equation}
[x,y,z]_C=\{x,y,z\}+\{y,z,x\}+\{z,x,y\},\quad \forall x,y,z\in A,\label{eq:3cc}
\end{equation}
\end{defi}

\begin{pro}{\rm \cite{baiguosheng}}
Let $(A,\{\cdot,\cdot,\cdot\})$ be a $3$-pre-Lie algebra. Then, $(A,[\cdot,\cdot,\cdot]_C)$ is a $3$-Lie algebra, where $[\cdot,\cdot,\cdot]_C$ is the induced $3$-commutator given by Eq.~\eqref{eq:3cc}. Furthermore, $L:\wedge^2A\longrightarrow\gl(A)$, which is defined by
\begin{equation}\label{eqn:repleft}
  L_{x,y}z=\{x,y,z\},\quad \forall~x,y,z\in A,
\end{equation}
is a representation of the $3$-Lie algebra $(A,[\cdot,\cdot,\cdot]_C)$ on $A$.
\end{pro}

\begin{defi}
Let $(\frkg,[\cdot,\cdot,\cdot] )$ be a $3$-Lie algebra. A {\bf symplectic structure} on $(\g,[\cdot,\cdot,\cdot])$ is a nondegenerate skew-symmetric bilinear form $\omega:\wedge^2\g\longrightarrow \mathbb R$, such that for all $x,y,z,w\in \g$, the following identity hold:
 \begin{equation}\label{eq:3sb}
\omega([x,y,z],w)-\omega([x,y,w],z)+\omega([x,z,w],y)-\omega([y,z,w],x)=0.
\end{equation}
\end{defi}

\begin{pro}{\rm \cite{baiguosheng}}\label{pro:inducedleft}
Let $(\frkg,[\cdot,\cdot,\cdot] )$ be a $3$-Lie algebra, and $\omega$  a  symplectic structure on $(\frkg,[\cdot,\cdot,\cdot] )$. Then, there exists
a compatible $3$-pre-Lie algebra structure on $\g$ given by
\begin{equation}\label{eq:bracketomega}
\omega(\{x,y,z\},w)=-\omega(z,[x,y,w]),\quad \forall x,y,z,w\in \g.
\end{equation}
\end{pro}

\begin{thm}
 Let $(\frkg,[\cdot,\cdot,\cdot] )$ be a $3$-Lie algebra,  $\omega$  a  symplectic structure on $(\frkg,[\cdot,\cdot,\cdot] )$, and
  $(\frkg,\{\cdot,\cdot,\cdot\})$  the induced $3$-pre-Lie
  algebra given in Proposition \ref{pro:inducedleft}. On the complex of
  vector spaces
  $\frkg^*\stackrel{{\omega^\sharp}^{-1}}{\longrightarrow}\frkg$,
  define $l_3$, which is totally skew-symmetric,  by
  \begin{equation}\label{eqn:homl2final}
\left\{\begin{array}{rlll}l_3(x,y,z)&=&[x,y,z], &\forall ~x,y,z\in \frkg\\
l_3(x,y,\xi)&=& L ^*_{ x,y}\xi, &\forall ~x,y\in \frkg, ~\xi\in \frkg^*\\
l_3(x,\xi,\eta)&=&0,&\forall ~x\in \frkg, ~\xi,\eta\in \frkg^*\\
l_3(\xi,\eta,\gamma)&=&0, &\forall~ \xi,\eta,\gamma\in \frkg^*,\end{array}\right.
\end{equation}
where $L ^* $ is the dual representation of
$L $ and $\omega^\sharp:\g\longrightarrow\g^*$ is given by $\omega^\sharp(x)(y)=\omega(x,y)$. Then,
$(\frkg^*,\frkg,\dM={\omega^\sharp}^{-1},l_3)$
is a strict $3$-Lie $2$-algebra.
\end{thm}
\pf Obviously, we only need to show that Conditions (a) and (c) in Lemma \ref{lem:2term 3L} hold,   the other conditions hold naturally. For all $x, y,z,w\in \g$,  let $f=\omega^\sharp(z)$ and $g=\omega^\sharp(w)$, we have
\begin{eqnarray*}
 \langle\dM l_3(x,y,f),g\rangle =-\langle l_3(x,y,\omega^\sharp(z)),w\rangle=-\langle L^*_{x,y}\omega^\sharp(z),w\rangle=\langle\omega^\sharp(z),\{x,y,w\}\rangle=\omega(z,\{x,y,w\}),
\end{eqnarray*}
and
\begin{eqnarray*}
 \langle l_3(x,y,\dM f),g\rangle =\langle l_3(x,y,z),\omega^\sharp(w)\rangle=  \omega(w,[x,y,z]).
\end{eqnarray*}
By \eqref{eq:bracketomega}, we deduce that $\dM l_3(x,y,f)=l_3(x,y,\dM f)$, i.e. Condition (a) hold.
Furthermore, also by \eqref{eq:bracketomega}, we have
\begin{eqnarray*}
 \langle l_3(\dM f,g,x),y\rangle&=&\langle L^*_{x,z}\omega^\sharp(w),y\rangle=-\omega(w,\{x,z,y\})=-\omega(y,[x,z,w]),\\
 \langle l_3( f,\dM g,x),y\rangle&=&\langle L^*_{w,x}\omega^\sharp(z),y\rangle=-\omega(z,\{w,x,y\})=-\omega(y,[w,x,z]).
\end{eqnarray*}
Therefore, we have $l_3(\dM f,g,x)=l_3( f,\dM g,x)$, which implies that Condition (b) hold. The proof is finished. \qed


\begin{thebibliography}{999}


\bibitem{ammardefiLeibnizalgebra}
 M. Ammar and N. Poncin, Coalgebraic Approach to the Loday Infinity Category, Stem
Differential for $2n$-ary Graded and Homotopy Algebras, \emph{Ann.
Inst. Fourier (Grenoble).} 60 (1) (2010),  355-387.



\bibitem{baez:2algebras}
 J. Baez and A. S. Crans, Higher-Dimensional Algebra VI: Lie
 2-Algebras, \emph{Theory and Appl.  Categ.} 12 (2004),
 492-528.

\bibitem{BL3}
J. Bagger and N. Lambert,  Three-algebras and N=6 Chern-Simons gauge theories. \emph{ Phys. Rev. D} 79 (2009), no. 2, 025002, 8 pp.

 \bibitem{baiguosheng}
 C. Bai, L. Guo and Y. Sheng, Bialgebras, classical Yang-Baxter equation and Manin triple for 3-Lie algebras, in preparation.

 \bibitem{BaiRGuo}
R. Bai, L. Guo, J. Li and Y. Wu, Rota-Baxter 3-Lie algebras, {\em J. Math. Phys.} {\bf 54} (2013), 064504, 14pp.


\bibitem{CasasPBW}
J. M. Casas, M. Insua and M. Ladra, Poincaré-Birkhoff-Witt theorem for Leibniz $n$-algebras. \emph{J. Symbolic Comput.} 42 (2007), no. 11-12, 1052-1065.


\bibitem{CasasSN}
J. M. Casas, E. Khmaladze, and M. Ladra, On solvability and nilpotency of Leibniz $n$-algebras. \emph{Comm. Algebra} 34 (2006), no. 8, 2769-2780.


\bibitem{CasasCMn-Lie}
J. M. Casas, E. Khmaladze, and M. Ladra, Crossed modules for Leibniz $n$-algebras. \emph{Forum Math.} 20 (2008), no. 5, 841-858.

\bibitem{CasasHigher}
J. M. Casas, E. Khmaladze, and M. Ladra, Higher Hopf formula for homology of Leibniz $n$-algebras. \emph{ J. Pure Appl. Algebra} 214 (2010), no. 6, 797-808.


\bibitem{BL4}
S. Cherkis and C. S{$\ddot{\rm a}$}mann,   Multiple M2-branes and generalized 3-Lie algebras. Phys. Rev. D 78 (2008), no. 6, 066019, 11 pp.

 \bibitem{review}
J. A. de Azcarraga and J. M. Izquierdo, $n$-ary algebras: a review with applications,
\emph{ J. Phys. A: Math. Theor.} 43 (2010), 293001.

\bibitem{Fialowski}
A. Fialowski and M. Penkava,  Deformation theory of infinity algebras. \emph{J. Algebra} 255 (2002), no. 1, 59-88.


\bibitem{deformation}
J. Figueroa-O$'$Farrill, Deformations of 3-algebras. \emph{J. Math. Phys.} 50 (2009), no. 11, 113514, 27 pp.

\bibitem{Filippov}  V. T. Filippov, $n$-Lie algebras,  {\it Sib. Mat. Zh.} 26 (1985) 126-140.

\bibitem{Gautheron}
P. Gautheron, Some remarks concerning Nambu mechanics, \emph{Lett. Math. Phys.}, 37 (1996) 103-116.


\bibitem{BL2}
J. Gomis, D. Rodríguez-Gómez, M. Van Raamsdonk and H. Verlinde,  Supersymmetric Yang-Mills theory from Lorentzian three-algebras. \emph{J. High Energy Phys.} 2008, no. 8, 094, 18 pp.

\bibitem{HHM} P. Ho, R. Hou and Y. Matsuo, Lie $3$-algebra and multiple
$M_2$-branes,  \emph{J. High Energy Phys.}  2008,  no. 6, 020, 30 pp.
\bibitem{Kasymov}
Sh. M. Kasymov, On a theory of n-Lie algebras. (Russian) \emph{Algebra i Logika} 26 (1987), no. 3, 277-297, 398.

 \bibitem{livernet}
M. Livernet, Homologie des alg$\rm\grave{e}$bres stables de matrices sur une $A_\infty$-alg$\rm\grave{e}$bre, \emph{C. R. Acad. Sci. Paris S$\rm\acute{e}$r. I
Math.}, 329(2) (1999), 113-116.



\bibitem{N} Y. Nambu, Generalized Hamiltonian dynamics, {\it Phys. Rev. D} 7 (1973)
                2405-2412.
                
 \bibitem{Penkava}                
M. Penkava,   L-infinity algebras and their cohomology,  arXiv:q-alg/9512014.
   


\bibitem{P} G. Papadopoulos, M2-branes, $3$-Lie algebras and
Plucker relations,  \emph{J. High Energy Phys.}  2008,  no. 5, 054, 9 pp.

\bibitem{Leibniz2al}Y. Sheng and Z. Liu, Leibniz $2$-algebras and twisted Courant algebroids.\emph{ Comm. Algebra } 41 (2013), no. 5, 1929–1953.

\bibitem{Stasheff1} M. Schlessinger and J. Stasheff, The Lie algebra
structure of tangent cohomology and deformation theory, \emph{J.
Pure Appl. Algebra} 38 (1985), 313-322.


\bibitem{T} L. Takhtajan, On foundation of the generalized Nambu mechanics,
{\it Comm. Math. Phys.} 160 (1994) 295-315.

\bibitem{UchinoshL}
K. Uchino, Derived brackets and sh Leibniz algebras, \emph{J. Pure
Appl. Algebra}, 215 (2011) 1102-1111.

\bibitem{BL1}
M. Van Raamsdonk,   Comments on the Bagger-Lambert theory and multiple M2-branes. \emph{J. High Energy Phys.} 2008, no. 5, 105, 9 pp.


\end{thebibliography}
\end{document}